\documentstyle[12pt]{article}
\begin{document}

\makeatletter
\renewcommand{\@biblabel}[1]{#1. \hfill}
\makeatother


\title{A new proof of the Frobenius conjecture
on the dimensions of real algebras
 without zero divisors}
\author{ K.E. Feldman}
\date{}
\maketitle

\noindent
\raisebox{5cm}[0pt][0pt]{Moscow Univ. Math. Bull.
55 (2000), no. 1, p. 48-50.}

A new way to prove the Frobenius conjecture on the dimensions of real
algebras without zero divisors is given in the present paper.
Firstly, the proof of nonexistence of real algebras without zero
divisors in all dimensions except 1,2,4 and 8 was given
in~\cite{Adams}. It was based on the simplicial cohomology operation
technique. Later on the methods of $K$-theory cohomology operations gave one
a possibility to obtain a more simple proof of the Frobenius conjecture
(see~\cite{Athiyah}).

For proving the Frobenius conjecture we suggest a new approach
different from~\cite{Adams,Athiyah}. We demonstrate that the
restriction on the dimensions of real algebras without zero divisors
follows elementary from the structure of $K$-functors of real projective
spaces. The general idea is to use $K$-theory characteristic classes
for investigation the question of parallelizibility of real projective
spaces that is equivalent to the Frobenius conjecture(see, e.g.
\cite{MilnorStashef}). Simplification of this scheme lies in the
basis of our proof. The structure of this paper is as follows. We begin with
the calculation of $K^*(RP^n,\emptyset)$. Then we give without proof the
reduction of the Frobenius conjecture to the question of parallelizibility
of real projective spaces. Finally, we obtain exact dimensions
of real algebras without zero divisors.

Let $\xi^1_n$ be the one-dimensional real Hopf vector bundle over $RP^n$,
and let $\xi^{\bot}_n$ be its orthogonal complement. Denote
the one-dimensional Hopf complex vector bundle over $CP^n$ by $\eta^1_n$.
We give the simplest calculation of $K^*(RP^n)$.

Our method to calculate $K$-functor of $RP^n$ is based on the
following geometric observation.
The complex Hopf bundle $\pi:S^{2n+1}\to CP^n$ with fiber $S^1$
can be passed through the real Hopf bundle $\pi_1:S^{2n+1}\to RP^{2n+1}$
with fiber $Z_2$. Under these conditions we obtain the bundle
$\pi_2:RP^{2n+1}\to CP^n$ whose fiber is also a circle.
More over, the following theorem holds.

{\bf Theorem 1.} {\it The bundle $\pi_2:RP^{2n+1}\to CP^n$ is isomorphic
to the spherical bundle of the tensor square of the bundle $\eta^1_n$
and also $\pi^*_2\eta^1_n\cong C\otimes \xi^1_{2n+1}$.}

{\bf Proof.} Denote the tensor square of the bundle $\eta^1_n$ by
$\eta^2$.  Let us construct an equivariant with respect to
$S^1$-action homeomorphism $g$ of the spaces $RP^{2n+1}$ and
$S(\eta^2)$.  Observe that $S(\eta^2)\cong S^{2n+1}\times_{\rho} S^1$,
where $\rho:S^1\times S^1\to S^1$ is defined by the formula
$\rho(u,v)=u^2v$. Construct the map $g:RP^{2n+1}\to S(\eta^2)$
supposing
$$
g(x_1:x_2:...:x_{2n+1}:x_{2n+2})=((\bar w z_1,...,\bar
w z_{n+1}),w^2)
$$
where $z_j=x_{2j-1}+ix_{2j}$, $j=1,...,n+1$, and $w$
runs over all complex numbers whose absolute values are 1. It is easy to see
that the map constructed is defined correctly, one--to--one,
continuous and equivariant with respect to $S^1$-action.
This completes the proof of the first part.

For proving the second part it is enough to observe
that the bundle $\eta^1_n$ can be defined as a set of proportions
$$
\{(z_1:...:z_{n+1}:\lambda)|\lambda,z_j\in
C,\sum^{n+1}_{j=1}|z_j|^2\ne 0\}.
$$
The lifting of this bundle on $RP^{2n+1}$ by the projection
$\pi_2$ is a set of the proportions
$$
\{(x_1:x_2:...:x_{2n+1}:x_{2n+2}:u:v)|z_j=x_{2j-1}+ix_{2j},
\lambda=u+iv, \sum^{2n+2}_{j=1}x^2_j\ne 0\}.
$$
The last bundle is exactly the bundle $C\otimes \xi^1_{2n+1}$.
This completes the proof of Theorem 1.

{\bf Theorem 2 \cite{Athiyah}.} {\it The group $K^0(RP^m,\emptyset)$
is isomorphic to the direct sum $Z\oplus Z_{2^{[m/2]}}$.
The generator of the second summand is given by the stable equivalence
class of the bundle $C\otimes \xi^1_m$. The group
$K^1(RP^m,\emptyset)$ equals to zero when $m$ is even, and
to $Z$ when $m$ is odd.}

 The proof of this theorem is based on the following well--known
lemma~\cite{Athiyah}.

{\bf Lemma.} {\it There is an isomorphism of groups
$K^0(CP^n,\emptyset)=Z[\beta]/\{\beta^{n+1}=0\}$, where
$-\bar\beta=[\eta^1_n]-[1]$ is the stable equivalence class
of the one-dimensional complex Hopf vector bundle $\eta^1_n$ over
$CP^n$. The group $K^1(CP^n,\emptyset)$ equals to zero.}

{\bf Remark.} In the proof of the lemma one may only use the Thom isomorphism
in complex $K$-theory and the fact that $CP^n$ is the Thom space of
the one-dimensional complex Hopf vector bundle over $CP^{n-1}$.

{\bf Proof of Theorem 2.} Firstly we investigate the case $m=2n+1$.
Consider the bundle constructed in Theorem 1: $\pi_2:RP^{2n+1}\to CP^n$.
That is a spherical bundle associated with the tensor square of the bundle
$\eta^1_n$ over $CP^n$. Let, as above, $\eta^2\cong\eta^1_n\otimes\eta^1_n$.
Consider the exact sequence of the pair:
$$
\begin{array}{ccccc}
K^0(D\eta^2,S\eta^2)    &\to  &K^0(CP^n,\emptyset)&\to  &K^0(RP^{2n+1},\emptyset)\\
\uparrow                &     &                   &     &\downarrow              \\
K^1(RP^{2n+1},\emptyset)&\gets&K^1(CP^n,\emptyset)&\gets&K^1(D\eta^2,S\eta^2).
\end{array}
\eqno(*)
$$

Note that the first and fourth terms of this sequence can be
changed on $K^0(T\eta^2)$ and $K^1(T\eta^2)$ respectively, where $T\eta^2$
is the Thom space of the bundle $\eta^2$. Besides, the image of the
homomorphism $K^0(T\eta^2)\to K^0(CP^n,\emptyset)$ coincides with
the image of the composition homomorphism
$$
K^0(CP^n,\emptyset)\cong K^0(T\eta^2)\to K^0(CP^n,\emptyset).
$$
The last homomorphism represents the multiplication
by the Euler class~\cite{Athiyah} of the bundle $\eta^2$, i.e. on
$[1]-[\bar\eta^2]=2[\beta]-[\beta]^2$.
As far as $K^1(T\eta^2)\cong K^1(CP^n,\emptyset)=0$,
we have $K^0(RP^{2n+1},\emptyset)\cong K^0(CP^n,\emptyset)/\{\beta^2=2\beta\}$.
Hence, $K^0(RP^{2n+1},\emptyset)\cong Z\oplus Z_{2^n}$.
Under these conditions the generator of the subgroup $Z_{2^n}$
can be given by element $-\pi^*_2(\beta)$. From Theorem 1
we obtain $-\pi^*_2(\beta)=[C\otimes\xi^1_{2n+1}]-[1]$.

Since $K^1(CP^n,\emptyset)=0$, from the exact sequence $(*)$
we have
$$
K^1(RP^{2n+1},\emptyset)\cong \mathrm{Ker}[K^0(T\eta^2)\to K^0(CP^n,\emptyset)].
$$
It is easy to see that
$\mathrm{Ker}[K^0(T\eta^2)\to K^0(CP^n,\emptyset)]\cong Z$.

 Owing to the functoriality of the exact sequence, the homomorphism
$ K^1(RP^{2n+1},\emptyset)\to
K^1(RP^{2n-1},\emptyset) $
is zero. Using the exact sequences of the pairs $(RP^{2n+1},RP^{2n})$ and
$(RP^{2n},RP^{2n-1})$ we find that the map
$K^1(RP^{2n+1},\emptyset)\to K^1(RP^{2n},\emptyset)$
is epimorphic while the map
$ K^1(RP^{2n},\emptyset)\to K^1(RP^{2n-1},\emptyset)$
is monomorphic. Thus, $K^1(RP^{2n},\emptyset)=0$. Then, from the
exact sequence of the pair $(RP^{2n+1},RP^{2n})$ it follows that
the map
$ K^0(RP^{2n+1},\emptyset)\to K^0(RP^{2n},\emptyset)$
is an isomorphism. This statement completes the proof.

The following theorem is the connecting bridge between Theorem 2 and the
Frobenius conjecture.

{\bf Theorem 3.} {\it Assume that there exists a bilinear operation of
multiplication $ p:R^n\times R^n\to R^n $ without zero divisors.
Then the bundle $\mathrm{Hom}(\xi^1_{n-1},\xi^{\bot}_{n-1}) $
is trivial.}

A simple proof of this theorem was suggested by Stiefel (see, for example,
\cite[\S 4]{MilnorStashef})

{\bf Theorem 4 \cite{Adams}.} {\it Real algebras without zero divisors
exist only in dimension $1,2,4$ and $8$.}

{\bf Proof.} Note that
$$
\mathrm{Hom}(\xi^1_{n-1},\xi^{\bot}_{n-1})\oplus 1\cong
\mathrm{Hom}(\xi^1_{n-1},\xi^{\bot}_{n-1}\oplus\xi^1_{n-1})\cong n\xi^1_{n-1}.
$$
Owing to Theorem 3, for existence of the algebra required in dimension $n$
it is necessary the bundle $ n\xi^1_{n-1} $ to be trivial. Consequently,
$n$ is divisible by the order of the element $[C\otimes \xi^1_{n-1}]-[1]$ in the group
$K^0(RP^{n-1},\emptyset)$.
According to Theorem 2, this means that $n$ is divisible by
$2^{[(n-1)/2]}$.
Let us write $n$ as $n=(2m+1)2^k$.  If $k=0$, then this condition is
only true for $m=0$, i.e. $n=1$. Let $k>0$. For $m>0$ we have
$[(n-1)/2]=(2m+1)2^{k-1}-1>k$, and the divisibility condition
can not be satisfied.
Thus, $m=0$. It is enough to observe that $2^{k-1}\le k+1$ only for
$k=1,2,3$. Hence, $n=1,2,4,8$.

It is well known that in dimension 1,2,4 and 8 there are real algebras without
zero divisors. One can choose, for example,  the field of real
numbers in dimension 1,  the field of complex numbers in dimension 2,
the quaternions in dimension 4 and  the Cayley algebra in dimension 8.

The author thanks L. A. Alaniya and V. M. Buchstaber for their helpful
comments and discussions.


Moscow State University

email: feldman@maths.ed.ac.uk

\end{document}